\documentclass[12pt]{elsarticle}
\usepackage[utf8]{inputenc}
\usepackage{natbib}

\usepackage{left_eq}   
\usepackage{algorithm}
\usepackage{algpseudocode}
\usepackage{etoolbox}

\algrenewcommand\algorithmicrequire{\textbf{Input:}}
\algrenewcommand\algorithmicensure{\textbf{Output:}}

\makeatletter
\newcommand*{\algrule}[1][\algorithmicindent]{%
  \makebox[#1][l]{%
    \hspace*{.2em}
    \vrule height .75\baselineskip depth .25\baselineskip
  }
}

\newcount\ALG@printindent@tempcnta
\def\ALG@printindent{%
    \ifnum \theALG@nested>0
    \ifx\ALG@text\ALG@x@notext
    \else
    \unskip
    \ALG@printindent@tempcnta=1
    \loop
    \algrule[\csname ALG@ind@\the\ALG@printindent@tempcnta\endcsname]%
    \advance \ALG@printindent@tempcnta 1
    \ifnum \ALG@printindent@tempcnta<\numexpr\theALG@nested+1\relax
    \repeat
    \fi
    \fi
}
\patchcmd{\ALG@doentity}{\noindent\hskip\ALG@tlm}{\ALG@printindent}{}{\errmessage{failed to patch}}
\patchcmd{\ALG@doentity}{\item[]\nointerlineskip}{}{}{} 
\makeatother

\usepackage{booktabs} 
\usepackage{multirow}
\usepackage{amsmath}
\usepackage{amssymb}
\usepackage{wrapfig}  
\usepackage{fancyhdr}
\usepackage{geometry}
\usepackage{graphicx} 
\usepackage{longtable}
\usepackage{hyperref}
\usepackage{setspace}

\hypersetup{
    colorlinks=false,
    linkcolor=blue,
    filecolor=magenta,      
    urlcolor=red,
    pdftitle={Overleaf Example},
    pdfpagemode=FullScreen,
    }

\makeatletter
\def\ps@pprintTitle{%
 \let\@oddhead\@empty
 \let\@evenhead\@empty
 \def\@oddfoot{}%
 \let\@evenfoot\@oddfoot}
\makeatother

\begin{document}
\onehalfspacing
\begin{frontmatter}

\title{Parametric Region Search: A Mixed-Integer Bilevel Optimization Problem Primal Heuristic}
\author[a]{Meng-Lin Tsai}
\author[a]{Parth Brahmbhatt}
\author[a]{Styliani Avraamidou\texorpdfstring{\corref{cor1}}{}}
\cortext[cor1]{Corresponding author Email: avraamidou@wisc.edu}

\affiliation[a]{Department of Chemical & Biological Engineering, University of Wisconsin-Madison, Madison, WI, 53706, USA}
\date{Feb 16 2026}

\begin{abstract}
  Bilevel optimization is a mathematical modeling formulation for hierarchical systems and two-player interactions, with wide-ranging applications in environmental, energy, and control engineering. Despite its utility, the mixed-integer bilevel optimization (MIBO) problem is exceptionally challenging to solve. While numerous exact and metaheuristic methods exist, the development of specialized primal heuristics for MIBO, aimed at quickly identifying high-quality feasible solutions, remains an underexplored area. This paper introduces the Parametric Region Search (PRS), a new primal heuristic for MIBO. The PRS method leverages insights from multi-parametric optimization by iteratively exploring regions defined by the lower-level problem's critical regions. We formally define the MIBO structure and the necessary parametric region formulations, and then detail the proposed heuristic's initialization and iterative search mechanism. Computational results demonstrate that the PRS heuristic consistently locates high-quality primal solutions compared to established derivative-free metaheuristics, including DOMINO-COBYLA and DOMINO-ISRES. Furthermore, we illustrate how the PRS can be effectively integrated with other heuristics like DOMINO-COBYLA to enhance the overall solution discovery process for MIBO.
\end{abstract}

\begin{keyword}
{Heuristic, Bilevel Optimization, Primal Heuristic, Multi-parametric optimization, Game Theory}
\end{keyword}
\end{frontmatter}

\section{Introduction}
As a mathematical framework for hierarchical interactions, bilevel optimization has become essential in fields ranging from industrial scheduling \citep{avraamidou2017multiparametric} and hierarchical control \citep{avraamidou2017multi} to electricity market pricing \citep{fampa2008bilevel}. However, the nested nature of these problems makes them NP-hard \citep{hansen1992new} and exceptionally difficult to solve. Solving bilevel problems is inherently complex, but the challenge intensifies significantly when the lower-level problem also contains discrete variables \citep{edmunds1992algorithm, kleinert2021survey}. To address these challenges, researchers have developed various approaches, ranging from decomposition techniques like Benders \citep{saharidis2009resolution}, and column generation \citep{yue2014game} to multiparametric (mp) optimization \citep{avraamidou2019multi, avraamidou2022multi,koppe2010parametric}, branch-and-bound\citep{moore1990mixed}, and branch-and-cut \citep{fischetti2018use,tahernejad2020branch}. Despite their theoretical rigor, these exact methods often fail to scale effectively when applied to large-scale problems. This computational bottleneck motivates a shift toward strategies that prioritize the rapid identification of high-quality feasible solutions rather than insisting on exact and global optimality. In this context, primal heuristics provide a practical alternative.

Primal heuristics are algorithms designed to quickly identify good feasible solutions for optimization problems, thereby establishing an upper bound on the optimal objective value. Although commonly utilized in standard Mixed-Integer Programming (MIP)\citep{berthold2006primal}, to our knowledge, the development and application of primal heuristics for the more complex domain of Mixed-Integer Bilevel Optimization (MIBO) are research areas that have received limited attention. 

Several heuristics have been developed to enhance the efficiency of solving bilevel optimization problems. For instance, the Improving Objective Cut heuristic \citep{denegre2011interdiction,tahernejad2020branch} identifies promising candidates by solving a restricted first-level problem, utilizing an objective cut from a known feasible follower response to narrow the search space. In contrast, the Second-level Priority heuristic \citep{denegre2011interdiction,tahernejad2020branch}, reverses this focus; rather than optimizing the first-level objective directly, it prioritizes the second-level objective to improve the chances of achieving bilevel feasibility. Beyond these approaches, local search techniques employing weak approximations of lower-level optimality have been proposed to facilitate the discovery of high-quality neighborhood solutions \citep{clementelocal2025}. Additionally, the DOMINO framework \citep{beykal2020domino} treats the lower-level optimization as a ``black box," allowing upper-level variables to be optimized using standard, off-the-shelf metaheuristic methods. The weighted sum heuristic \citep{tahernejad2020branch, denegre2011interdiction} combines the upper- and lower-level objectives using a variable weight to relax the bilevel optimization problem into a bi-objective optimization; it then solves for the lower-level responses to obtain a bilevel-feasible primal solution. Similarly, Goal Programming for Bi-level Optimization (GPBLO) \citep{shams2023designing} introduces a normalization step prior to weighting. GPBLO produces primal solutions that outperform Particle Swarm Optimization (PSO) in terms of both solution quality and computational efficiency in the computational study.

This paper introduces Parametric Region Search (PRS), a new heuristic designed to solve Mixed-Integer Bilevel Optimization (MIBO) problems with binary variables at the upper and lower level. The PRS approach is inspired by the multiparametric reformulation \citep{avraamidou2019multi2, avraamidou2019multi, avraamidou2022multi, koppe2010parametric} The upper-level parameter space is partitioned into critical regions (CRs), each defined by a set of linear inequalities that form a convex polyhedral set. Within each region, the set of active constraints remains constant. This allows the lower-level optimal value to be expressed as a piecewise linear function of the upper-level variables by the sensitivity theorem. By substituting the lower-level optimal value with the upper-level variable in each critical region, the original bi-level problem is transformed into a series of single-level optimization problems. A primary challenge with this method is that the number of critical regions grows exponentially as the problem size increases \citep{saini2025iteration}. This rapid expansion creates a bottleneck, resulting in poor scalability when attempting to solve larger problems.

Instead of exhaustively enumerating all critical regions in the entire parameter space as required by multiparametric reformulation methods, our proposed PRS approach selectively targets individual regions likely to contain high-quality bilevel feasible solutions. This strategy addresses a key inefficiency in standard multiparametric programming, where many identified critical regions fail to provide competitive solutions yet still require significant computational effort to resolve. By focusing the search on restricted high-potential single-level solution spaces while leveraging the underlying structure provided by multiparametric programming, the PRS avoids the overhead of a full reformulation. This approach scales more effectively than exhaustive methods, which often become computationally prohibitive. This targeted strategy reduces the required processing time, though it trades the guarantee of global optimality for increased practical efficiency.

The remainder of the paper is structured as follows. Section 2 formally defines the mixed-integer bilevel problem structure and the formulation of the multiparametric critical regions derived from the lower-level problem. Section 3 then describes the PRS heuristic, outlining its initialization and iterative region-based search mechanism. In Section 4, computational results are presented, comparing the PRS method against established local and global derivative-free metaheuristics. This comparison demonstrates that the PRS approach, utilizing fundamentally different logic rooted in structural analysis, consistently locates high-quality solutions. Section 5 explores two use cases, illustrating how the PRS can be coupled with other heuristics to enhance the discovery of better primal solutions for challenging mixed-integer optimization problems. Finally, Section 6 concludes the article with a discussion of the method's contribution and outlines promising avenues for future research in bilevel optimization.

\section{Bilevel Optimization and Multi-parametric (mp) Optimization}
This section establishes the theoretical context needed for the proposed PRS heuristic by introducing bilevel optimization and the multiparametric reformulation \citep{avraamidou2019multi,koppe2010parametric} that motivates this work.
\subsection{MILP-MILP bilevel optimization formulation}
The general Mixed-Integer Linear Program MILP-MILP bilevel problem can be formulated as follows:
\begin{equation}
\begin{aligned}
    & \min_{x_C, x_I, y_C, y_I} \quad F(x, y) = c_1^\top x + d_1^\top y \\
    & \text{s.t.}\quad A_1 x + B_1 y \leq b_1 \\
    & \quad \quad y \in \arg \min_{z_C, z_I} \quad f(x, z) = c_2^\top x + d_2^\top z \\
    & \quad \quad \text{s.t.} \quad A_2 x + B_2 z \leq b_2 \\
    & \text{where} \quad x_C \in \mathbb{R}^{n_{xC}}, \ x_I \in \mathbb{Z}^{n_{xI}} \\
    & \quad \quad \quad y_C,z_C \in \mathbb{R}^{n_{yC}}, \ y_I,z_I \in \mathbb{Z}^{n_{yI}}
\end{aligned}
\end{equation}

In this formulation, the upper-level variables are denoted by $x = (x_C, x_I)$, where $x_C \in \mathbb{R}^{n_{xC}}$ represents continuous decisions and $x_I \in \mathbb{Z}^{n_{xI}}$ represents discrete decisions. $z$ is the dummy variable for the lower-level problem to avoid a circular definition of the lower-level variable $y$. Both are similarly partitioned into continuous ($y_C,z_C$) and integer ($y_I,z_I$) components. The upper- and lower-level constraints, $A_1 x + B_1 y \leq b_1$ and $A_2 x + B_2 z \leq b_2$, define the feasible region. While $F(x, y) = c_1^\top x + d_1^\top y$ and $f(x, z) = c_2^\top x + d_2^\top z$ define the upper- and lower-objective. Note that equality constraints are omitted for brevity, as they can be expressed as pairs of inequalities.
The proposed framework requires all lower-level integer variables, $z_I$, to be binary. If the lower-level variables are general integers, they can be converted into binary form using standard expansion techniques, such as adding binary slack variables. Consequently, the terms ``integer'' and ``binary'' are used interchangeably throughout the remainder of this text.

\subsection{Lower-level subproblem (Inner Problem)}

For a fixed upper-level decision $\hat{x} = (\hat{x_C}, \hat{x_I})$, the lower-level problem reduces to a standard MILP, with the $c_2^\top \hat{x}$ term in the objective function being a constant that can be neglected. 

To analyze this subproblem, we decompose the lower-level cost vector $d_2$ and the constraint matrix $B_2$ based on the variable types. Specifically, we let $d_2 = [d_{2C}^\top, d_{2I}^\top]^\top$ and $B_2 = [B_{2C}, B_{2I}]$, where the subscripts $C$ and $I$ correspond to the continuous and integer (binary) dimensions, respectively. The inner problem is then expressed as:

\begin{equation}
\begin{aligned}
    & \min_{y_C, y_I} \quad d_{2C}^\top y_C + d_{2I}^\top y_I \\
    & \text{s.t.} \quad B_{2C} y_C + B_{2I} y_I \leq b_2 - A_{2} \hat{x} \\
    & \quad \quad y_C \in \mathbb{R}^{n_{yC}}, \ y_I \in \mathbb{Z}^{n_{yI}}
\end{aligned}
\label{eq:LLS}
\end{equation}
In Equation \eqref{eq:LLS}, the term $b_2 - A_{2} \hat{x}$ represents the effective right-hand side, showing how the leader's decision $\hat{x}$ restricts the follower's feasible space.

A bilevel problem is only feasible if this lower-level subproblem yields an optimal solution. However, if the lower level is degenerate, meaning multiple solutions $y$ minimize the follower's objective but result in different values for $F(x, y)$, the tie-breaking rule becomes critical. To handle this non-uniqueness, two standard approaches are used: the optimistic approach, which selects the variables that minimize the upper-level objective, and the pessimistic approach, which selects those that maximize it (the worst-case scenario) \citep{wiesemann2013pessimistic}. 

\subsection{Multiparametric (mp) Optimization Overview}
\label{Sec. mp_CR}

Before detailing the specific application to bilevel problems, we first introduce the general framework of multiparametric programming. The optimization formulation of a general mp linear programming (mp-LP) problem can be expressed as:
\begin{subequations}
\begin{align}
Z(\theta) = &\min_{t} \ c^\top t + \theta^\top H^\top t \\
\text{s.t.} &~At \leq b + F\theta \\
&\theta \in \Theta := \{ \theta \in \mathbb{R}^q \mid A_{\theta}\theta \leq b_{\theta} \} \\
&t \in \mathbb{R}^{n_c} \times \mathbb{Z}^{n_i}
\end{align}
\label{eq: mpLPproblem}
\end{subequations}
\noindent Here, the objective function is linear in the decision variables $t$ and is influenced by the parameters $\theta$ through the term $\theta^\top H^\top t$. The inequality constraints $At \leq b + F\theta$ show that the right-hand side limits are affine functions of the parameters $\theta$. The constraint $\theta \in \Theta$ defines the space of feasible parameters as a polytope bounded by the linear inequalities $A_{\theta}\theta \leq b_{\theta}$. In this general formulation, the decision vector $t$ can contain both continuous and integer variables, allowing the problem to be structured as a MILP. This formulation seeks an optimal solution $t$ for every possible value of $\theta$ within its defined bounds $\Theta$.

For instances where the problem is continuous (or once the integer variables in $x$ are fixed), the explicit parametric solution to Eq.~\eqref{eq: mpLPproblem} maps the parameter space into distinct geometric partitions known as critical regions. The solution can be expressed as a piecewise affine function of the parameters:

\begin{equation}
 \begin{aligned}
     {t}^*(\theta) = A^v\theta + b^v \quad \text{if} \quad \theta \in CR^v := \{ \theta \mid E^v\theta \leq f^v \} \\
     \forall \ v = 1, 2, 3, \dots, n_{CR}
 \end{aligned}
 \label{eq: mpSol}
 \end{equation}
\noindent where $t^{*}(\theta)$ is the optimal solution profile, $A^v$ and $b^v$ define the affine function for a specific critical region $v$, $E^v$ and $f^v$ represent the inequality matrices defining the critical region polytope $CR^v$, and $n_{CR}$ denotes the total number of critical regions. 

\subsection{Critical Region Formulation}
In the context of the bilevel optimization framework addressed in this work, the lower-level subproblem operates exactly as a multiparametric programming problem. Specifically, the upper-level decisions $x$ act as the varying parameters $\theta$, while the lower-level variables $y$ act as the decision variables.

Once the lower-level integer variables $y_I$ are fixed to values $\hat{y}_I$, and a feasible upper-level point $x_0$ is selected, the lower-level problem reduces to a standard linear programming (LP) problem. We can reformulate the bilevel problem into a number of single-level problems in the associated CRs using the multiparametric linear programming. 

\subsubsection{Identification of the Active Set}

Consider the fixed upper-level nominal point $x_0$, which transforms the multiparametric problem into a standard LP. Solving this LP yields the nominal lower-level solution $y_{C,0}$ and the corresponding Lagrange multipliers $\lambda_0$. Based on these values, the Active Set $\mathcal{A}$ is identified as the set of indices for constraints that satisfy strict equality:

\begin{equation}
\begin{aligned}
    B_{2C, \mathcal{A}} y_{C,0} = \hat{b}_{\mathcal{A}} - A_{2, \mathcal{A}} x_0 \\
    \lambda_{0, \mathcal{A}} > 0
\end{aligned}
\end{equation}
where $\hat{b}_{2,\mathcal{A}} = b_{2, \mathcal{A}} - B_{2I, \mathcal{A}} \hat{y}_I$ represents the effective right-hand side accounting for the fixed lower-level integers. We assume the Linear Independence Constraint Qualification (LICQ) holds, meaning the matrix $A_{2, \mathcal{A}}$ has full rank and is invertible.

\subsubsection{Parametric Solution}

According to the Basic Sensitivity Theorem for linear programming, the optimal basis (defined by $\mathcal{A}$) remains invariant in a neighborhood around the nominal point $x_0$. Consequently, the optimal continuous variables $y_C(x)$ can be expressed as an affine function of the parameter $x$:

\begin{equation}
y_C^*(x) = B_{2C, \mathcal{A}}^{-1} \left( \hat{b}_{2,\mathcal{A}} - A_{2, \mathcal{A}} x \right)
\end{equation}

Furthermore, for multiparametric LP (mp-LP) problems, the Lagrange multipliers $\lambda(x)$ do not change as a function of the parameter $x$ within the critical region; they remain equal to the nominal values $\lambda_0$.

\subsubsection{Generation of the Critical Region (CR)}

To ensure the affine function $y_C^*(x)$ remains the optimal solution, the solution must remain feasible. Since the dual variables $\lambda(x)$ are constant and non-negative by definition at the nominal point, the dual feasibility condition is naturally satisfied. Therefore, the critical region is defined strictly by the Primal Feasibility requirements of the inactive constraints. Substituting the parametric solution $y_C^*(x)$ into the original constraints yields the polyhedral representation of the critical region:

\begin{equation}
CR_{\mathcal{A}} = \left\{ x \in \mathbb{R}^{n_x} \ \middle| 
\begin{array}{l}
\ \left[ -B_{2C} B_{2C, \mathcal{A}}^{-1} A_{2C, \mathcal{A}} + A_{2} \right] x \leq \left[ \hat{b_2}_{\mathcal{A}} - B_{2C} B_{2C, \mathcal{A}}^{-1} \hat{b}_{\mathcal{A}} \right] \\
\ \left[ -B_{1C} B_{2C, \mathcal{A}}^{-1} A_{2C, \mathcal{A}} + A_{1} \right] x \leq \left[ \hat{b_1}_{\mathcal{A}} - B_{1C} A_{2C, \mathcal{A}}^{-1} \hat{b}_{\mathcal{A}} \right]
\end{array}
\right\}
\end{equation}

This formulation is modified from Eq. (2.7) in \citep{pistikopoulos2020multi}, adapted to the bilevel optimization context such that the parameter $\theta$ corresponds to upper-level variable $x$ and the constraint matrices are mapped to the bilevel problem structure. Within this compact polytope, the lower-level reaction $y_C^*(x)$ is uniquely determined and linear, allowing for the efficient single-level reformulation of the bilevel problem.

\subsubsection{mp-Reformulated Single-Level Problem}
Within a Critical Region $CR$, the optimal continuous response $y_C^*(x)$ is affine in $x$:
\begin{equation}
y_C^*(x) = B_{2C, \mathcal{A}}^{-1} \left( \hat{b}_{2,\mathcal{A}} - A_{2, \mathcal{A}} x \right) = K_C x + h_C
\end{equation}

Substituting $y_C^*(x)$ and the assumed fixed integer solution $\hat{y}_I$ into the upper-level problem yields the single-level mp-reformulated MILP:
\begin{equation}
\begin{aligned}
    & \min_{x_C, x_I} \quad F(x, y^*(x)) \\
    & \text{s.t.} \quad A_1 x + B_{1C} (K_C x + h_C) + B_{1I} \hat{y}_I \leq b_1 \\
    & \quad \quad x \in CR \\
    & \quad \quad x_C \in \mathbb{R}^{n_{xC}}, \ x_I \in \mathbb{Z}^{n_{xI}}
\end{aligned}
\end{equation}

The final bilevel solution is determined by identifying the minimum upper-level objective value across all critical regions. However, the number of these regions can grow exponentially with the bilevel optimization problem size \citep{saini2025iteration}, making the exhaustive enumeration of all CRs computationally prohibitive for large-scale applications. This challenge motivates the present work, which introduces a novel heuristic designed to strategically explore promising critical regions rather than enumerating them all, efficiently identifying high-quality solutions.

\section{Method}
\label{sec:Parametric Region Search Heuristic}
The Parametric Region Search (PRS) is a new primal heuristic developed to efficiently identify high-quality feasible solutions for MIBO problems given an initial starting point. To obtain a good starting point, Subsection \ref{Sec: Initialization} details the high-point relaxation response method. The heuristic is a three-step process described in Subsection \ref{Sec: PRS}, which consists of lower-level optimization, critical region (CR) generation, and regional upper-level optimization.

\subsection{Initialization (High-point Relaxation Response)}
\label{Sec: Initialization}
To obtain a starting feasible point, we solve the High-Point Relaxation (HPR) and determine the follower's response according to the procedure in Algorithm \ref{alg:hprr}. HPR is a bilevel optimization relaxation that removes the lower-level optimality constraints. By solving the resulting single-level MILP, a divided solution pair $(\hat{x}, \hat{y})$ is identified within the relaxed feasible region $\Omega$. Since the relaxed response $\hat{y}$ may not be optimal for the follower at the fixed leader decision $\hat{x}$, a feasibility response step is executed. The leader's decision is fixed, and the lower-level problem is solved in isolation to determine the true optimal response $y^*$. This process yields a bilevel feasible point $(\hat{x}, y^*)$, providing a starting point for local or global search heuristics. While Step 1 is sufficient for determining the initial upper-level variables $\hat{x}$, the subsequent Step 2 is included as an optional procedure to facilitate the gap analysis described in Section \ref{Sec: Quality}.

\begin{algorithm}
\caption{High-Point Relaxation (HPR) Response Initialization}
\label{alg:hprr}
\begin{algorithmic}[1]
\Require Bilevel formulation parameters $(c_1, d_1, A_1, B_1, b_1, c_2, d_2, A_2, B_2, b_2)$
\Ensure Initial feasible candidate $(\hat{x}, y^*)$

\State \textbf{Step 1: Relax Lower-Level Optimality}
\State Solve the relaxed upper-level MILP as a single-level problem:
$$ (\hat{x}, \hat{y}) \leftarrow \arg\min_{x,y} \{ c_1^\top x + d_1^\top y \mid (x,y) \in \Omega \} $$
\Comment{Where $\Omega$ represents all upper and lower constraints}

\State \textbf{Step 2: Feasibility Response (Optional)}
\State Fix $x \leftarrow \hat{x}$ and solve the lower-level problem:
$$ y^* \leftarrow \arg\min_{y} \{ f(\hat{x}, y) \mid g(\hat{x}, y) \leq 0 \} $$

\State \Return Initial candidate $(\hat{x}, y^*)$
\end{algorithmic}
\end{algorithm}

\subsection{Parametric Region Search (PRS) Heuristic}
\label{Sec: PRS}
\begin{figure}[htbp]
    \centering
    \includegraphics[width=0.95\linewidth]{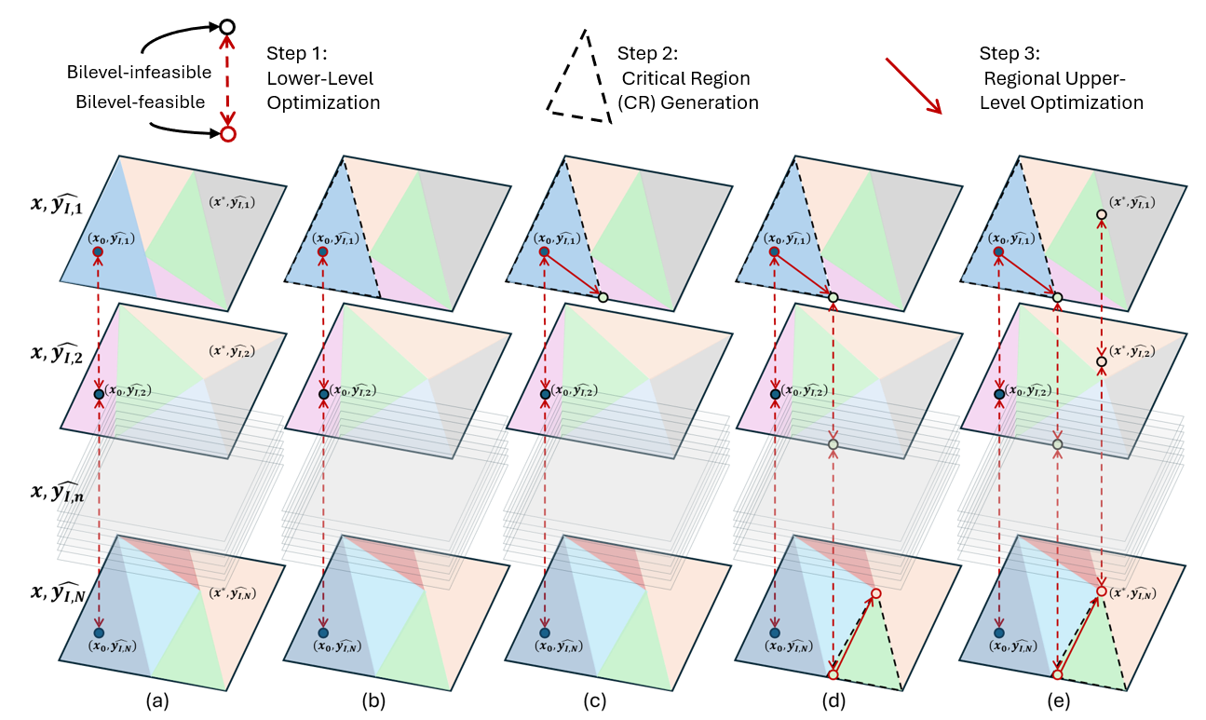}
    \caption{Schematic diagram of PRS heuristic. (a) Step 1, lower level optimization with a fixed set of upper variables, highlighted in red dashed vertical arrows. (b) Step 2, Critical region generation, highlighting in dashed border lines. (c) Step 3, regional upper level optimization, highlighted by the red arrow. (d) Repeat Step 1, Step 2, and Step 3 in a new iteration. (e) Repeat Step 1 and converge.}
    \label{fig:Schematic diagram}
\end{figure}

This subsection presents the proposed methodology, which is detailed in Figure \ref{fig:Schematic diagram}. Following the methodology and providing a detailed step-by-step procedure via the pseudo-code in Algorithm \ref{alg:prs}. We conclude with a numerical example to illustrate the process. The methodology utilizing a multi-parametric reformulation structure involves $N$ variants of the upper-level variable $x$-space, where $N$ represents the total number of possible combinations of lower-level binary variables. Within the schematic representation, these variants are visualized as stacked planes, where each plane $(x, \widehat{y_{I,n}})$ denotes an upper-variable space corresponding to a specific configuration of lower-level integer variable combinations. Each variant is further partitioned into multiple critical regions, represented by different colors in Figure \ref{fig:Schematic diagram}. A challenge of solving this structure is that $N$ grows combinatorially with the number of lower-level binary variables, which can lead to an exceptionally large search space. While traditional multi-parametric bilevel solvers typically require an exhaustive search of every critical region across all planes to ensure global optimality \citep{avraamidou2019multi, koppe2010parametric}, the algorithm presented here iterates only through selected regions. 

The algorithm starts with a fixed upper-level point $x_{0}$, which is established using the high-point relaxation response, algorithm \ref{alg:hprr}, discussed in Section \ref{Sec: Initialization}. To verify bilevel feasibility, Step 1.1 involves solving the lower-level subproblem (defined in Eq. \ref{eq:LLS}) to identify the true follower response. In Figure \ref{fig:Schematic diagram}(a), this lower-level optimization process is represented by vertical dashed red lines that span the stacked planes. The procedure identifies a bilevel feasible point, which is distinguished by a red-bordered circle. In contrast, other points sharing the same fixed value $x$ but possessing different $\hat{y}_{I,n}$ values are identified as bilevel infeasible and are marked with black borders.

In Step 1.2, the algorithm performs a lexicographic update to maintain the current best solution. The candidate solution $(x_{k}, y_{k})$ is accepted as the new global optimum $(x^{*}, y^{*})$ if it provides a strictly better upper-level objective value ($F_{curr} < F_{best}$). In the event of a tie in the upper-level objective (degenerate solution in upper-level), the lower-level objective value is used as a tie-breaker, where the solution is updated if $f_{curr} < f_{best}$. Subsequently, Step 1.3 executes a convergence check to determine if the search should terminate. This is implemented if the newly identified point falls within a previously visited critical region, at which point the local optimization phase terminates.

In Step 2, the algorithm focuses on critical region generation by utilizing the bilevel feasible points previously identified in Step 1. Following the methodology detailed in Section \ref{Sec. mp_CR}, this process constructs a polyhedral region, denoted as $CR_{\mathcal{A}}$, which is represented by dashed boundaries in Figure \ref{fig:Schematic diagram}(b). Within this specific region, the follower's optimal active constraints set $\mathcal{A}_k$ remains invariant, allowing the lower-level response to be characterized as an explicit parametric function $y_C(x)$. Consequently, the complex bilevel problem is reduced to a more manageable single-level program that is restricted to the current polyhedral subset. 

In Step 3, the algorithm executes a regional single-level optimization by performing a local search within the current critical region. This step facilitates the transition using regional upper-level optimization (represented by the red arrow) from a bilevel feasible point ($x_{k}$; represented in the schematic by a red-outlined circle) toward a potential candidate optimum ($x_{new}$; represented by a hollow circle), as illustrated in Figure \ref{fig:Schematic diagram}c. It is important to note that this new candidate point may not necessarily be a bilevel feasible solution. This occurs because the critical region guarantees lower-level optimality only with respect to a fixed lower-level active constraint set $\hat{y}_{I,n}$, rather than the global problem. 

Following this local search, the algorithm returns to Step 1 to re-evaluate the true follower response for the newly identified $x_{new}$, as illustrated in Figure \ref{fig:Schematic diagram}d. In Figure \ref{fig:Schematic diagram}e, the algorithm reaches termination as the identified new point resides within a previously solved CR, signifying that no further exploration of this area is required.

Figure \ref{fig:Schematic diagram} illustrates a fundamental distinction between the proposed Parametric Regional Search (PRS) approach and traditional multiparametric reformulation methods. While each plane in the decision space is characterized by distinct critical regions, the PRS approach does not aim to exhaustively identify all critical regions in the manner of the mp-reformulation method \citep{avraamidou2019multi,koppe2010parametric}. Instead, the algorithm functions more selectively by locating only the specific region containing the identified bilevel-feasible solution during each iteration. The PRS is guaranteed to terminate because the number of CRs is finite. The bilevel feasible solutions identified by the PRS are for the optimistic formulation of bilevel optimization. This is related to the search mechanism, which prioritizes improvements in the upper-level objective when exploring new critical regions in Step 1.

\begin{algorithm}
\caption{Parametric Region Search (PRS)}
\label{alg:prs}
\begin{algorithmic}[1]
\Require Initial point $x_0$ from Algorithm \ref{alg:hprr}, Max iterations $k_{max}$
\Ensure Best found solution $(x^*, y^*)$

\State $k \leftarrow 0$
\State $F_{best} \leftarrow \infty, f_{best} \leftarrow \infty$ \Comment{Initialize best upper and lower objectives}
\State $x_k \leftarrow x_0$

\While{$k < k_{max}$}
    \State \textbf{Step 1.1: Lower-Level Optimization}
    \State Solve $y_k \in \arg\min_y \{ f(x_k, y) \mid (x_k, y) \in \Omega \}$
    \State $F_{curr} \leftarrow F(x_k, y_k)$
    \State $f_{curr} \leftarrow f(x_k, y_k)$

    \State \textbf{Step 1.2: Update Best Solution (Lexicographic)}
    \If{$F_{curr} < F_{best}$ \textbf{or} ($F_{curr} = F_{best}$ \textbf{and} $f_{curr} < f_{best}$)}
        \State $(x^*, y^*) \leftarrow (x_k, y_k)$
        \State $F_{best} \leftarrow F_{curr}, f_{best} \leftarrow f_{curr}$
    \EndIf

    \State \textbf{Step 1.3: Convergence Check}
    \If{$\exists i \in \{0, \dots, k-1\} : x_k \in CR_i$} \Comment{Revisiting explored CRs}
    \State \textbf{break}
    \EndIf
    
    \State \textbf{Step 2.1: Critical Region (CR) Generation}
    \State Perform sensitivity analysis (Sec. \ref{Sec. mp_CR}) at $x_k$ to identify the active set $\mathcal{A}_k$.
    \State Construct $CR_k = \{ x \mid \mathcal{A}_k \text{ remains optimal for the lower level} \}$.

    \State \textbf{Step 3: Regional Upper-Level Optimization}
    \State Set $y_C(x)=K_C x + h_C$ for $x \in CR_k$.
    \State Solve: $x_{new} \leftarrow \arg\min_{x} \{ F(x, y(x)) \mid x \in CR_k \}$
    

    \State $x_{k+1} \leftarrow x_{new}$
    \State $k \leftarrow k + 1$
\EndWhile
\State \Return $(x^*, y^*)$
\end{algorithmic}
\end{algorithm}
\subsection{Numerical example:}

\begin{figure}[htbp]
    \centering
    \includegraphics[width=0.5\linewidth]{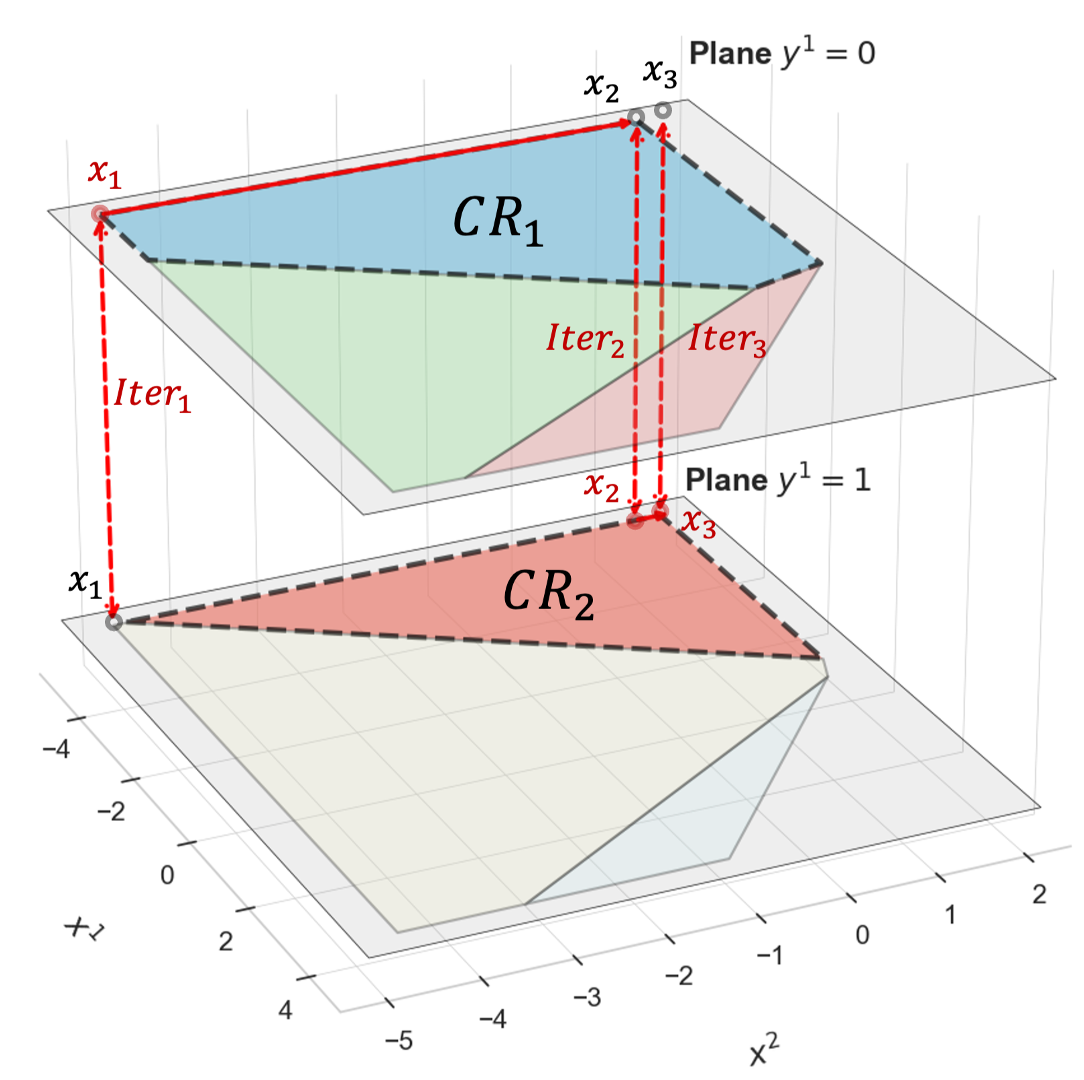}
    \caption{Schematic diagram for the numerical example.}
    \label{fig:numerical eg}
\end{figure}

The LP-MILP numerical example (\ref{eq:num_example}) is solved bellow to illustrate the steps of the proposed algorithm. Fig. \ref{fig:numerical eg} is schematically presenting the iterations of the algorithm.

\begin{equation}
\label{eq:num_example}
\begin{aligned}
\min_{x, y} \quad & F(x, y) = \begin{bmatrix} 0 & 28.1 \end{bmatrix} x + \begin{bmatrix} 0 & 26.2 \end{bmatrix} y \\
\text{s.t.} \quad & \begin{bmatrix} 1 & 0 \\ 0 & 1 \\ -1 & 0 \\ 0 & -1 \end{bmatrix} x \le \begin{bmatrix} 4.85 \\ 4.85 \\ 4.85 \\ 4.85 \end{bmatrix} \\
& y \in \arg \min_{z} f(x, z) = \begin{bmatrix} 1 & -3.6 \end{bmatrix} z \\
& \text{s.t.} \quad \begin{bmatrix} -2.0 & 5.9 \\ -3.0 & -3.0 \\ 7.1 & 4.7 \\ 0 & 2.0 \\ 0 & -10.0 \end{bmatrix} z \le \begin{bmatrix} 3 \\ 14.6 \\ 6.2 \\ 10.7 \\ 100 \end{bmatrix} + \begin{bmatrix} 0 & -7.8 \\ 0 & -9 \\ -4.7 & -2.4 \\ -7.8 & -5.9 \\ 0 & 0 \end{bmatrix} x \\
\text{where} \quad & x \in \mathbb{R}^2 \\
& y^1, z^1 \in \{ 0, 1 \}, \quad y^2, z^2 \in \mathbb{R}
\end{aligned}
\end{equation}

Here, we introduce the notation: $V_{k}^{index}$. In this notation, the superscript represents the index within the $x$ or $y$ vector, while the subscript, $k$, denotes the iteration number.

\subsubsection{Initialization}

An initial guess $x_0$ is required as an input to Algorithm \ref{alg:prs}. Algorithm \ref{alg:hprr} (HPR) is followed to generate this initial guess, where the bilevel problem is treated as a single-level MILP by ignoring the lower-level optimality constraint, resulting to the MILP formulation (Eq. \ref{eq:num_example_HPR}). Note that the variable $z$ is used instead of $y$ to indicate that the lower-level variable in this subproblem is not solved to lower-level optimality.

\begin{equation}
\label{eq:num_example_HPR}
\begin{aligned}
\min_{x_0, z} \quad & F(x_0, z) = \begin{bmatrix} 0 & 28.1 \end{bmatrix} x_0 + \begin{bmatrix} 0 & 26.2 \end{bmatrix} z \\
\text{s.t.} \quad & \begin{bmatrix} 1 & 0 \\ 0 & 1 \\ -1 & 0 \\ 0 & -1 \end{bmatrix} x_0 \le \begin{bmatrix} 4.85 \\ 4.85 \\ 4.85 \\ 4.85 \end{bmatrix} \\
&\begin{bmatrix} -2.0 & 5.9 \\ -3.0 & -3.0 \\ 7.1 & 4.7 \\ 0 & 2.0 \\ 0 & -10.0 \end{bmatrix} z \le \begin{bmatrix} 3 \\ 14.6 \\ 6.2 \\ 10.7 \\ 100 \end{bmatrix} + \begin{bmatrix} 0 & -7.8 \\ 0 & -9 \\ -4.7 & -2.4 \\ -7.8 & -5.9 \\ 0 & 0 \end{bmatrix} x_0 \\
\text{where} \quad & x_0 \in \mathbb{R}^2 \\
& z^1 \in \{ 0, 1 \}, \quad z^2 \in \mathbb{R}
\end{aligned}
\end{equation}

Solving \ref{eq:num_example_HPR}, the HPR yields an initial upper-level vector $x_0 = [-4.85 -4.85]$. This point serves as the input for the first iteration in Algorithm \ref{alg:hprr} (PRS). 

\subsubsection{Iteration 1}

The procedure starts at an initial parameter $x_0$. In step 1, the optimization of the lower-level MILP problem Eq. \ref{eq:Iter1_Step1} yields the optimal response $y_1 = [0,6.9203]$. 
\begin{equation}
\label{eq:Iter1_Step1}
\begin{aligned}
\min_{y_1} \quad & f(y_1) = \begin{bmatrix} 1 & -3.6 \end{bmatrix} y_1\\
\text{s.t.} \quad & \begin{bmatrix} -2.0 & 5.9 \\ -3.0 & -3.0 \\ 7.1 & 4.7 \\ 0 & 2.0 \\ 0 & -10.0 \end{bmatrix} y_1 \le \begin{bmatrix} 3 \\ 14.6 \\ 6.2 \\ 10.7 \\ 100 \end{bmatrix} + \begin{bmatrix} 0 & -7.8 \\ 0 & -9 \\ -4.7 & -2.4 \\ -7.8 & -5.9 \\ 0 & 0 \end{bmatrix} x_0\\
\text{where} \quad & y^1_1 \in \{ 0, 1 \}, \quad y^2_1 \in \mathbb{R}
\end{aligned}
\end{equation}
Since this represents the initial bilevel feasible solution, now $x_0 = x_1$ and the pair $(x_1, y_1)$ is stored as the incumbent best solution $(x^*, y^*)$, with an upper objective of $45.03$. In step 2, a sensitivity analysis is performed at $x_1; y_1^1 = 0$ to identify the active set $\mathcal{A}_1 = \{0\}; y^1=0$, signifying that the zeroth lower-level constraint is active and the lower-level binary is fixed at 1. Within the resulting polyhedral region $CR_1$, defined by $E_1x \leq f_1$, and the lower-level continuous variable, it is expressed as a linear map $y^2(x) = K_{C1} x + h_{C1}$. Where:$$K_{C1} = \begin{bmatrix} 0 & -1.3220 \end{bmatrix}, \quad h_{C1} = \begin{bmatrix} 0.5085 \end{bmatrix}$$
The polyhedral constraints defining the $CR_1$ are defined by $E_1x \leq f_1$, where:$$E_1 = \begin{bmatrix} 0 & 1 \\ 0.7765 & -0.6301 \\ 0.9228 & 0.3852 \\ -1 & 0 \\ 0 & -1 \end{bmatrix}, \quad f_1 = \begin{bmatrix} 1.2437 \\ 0.6295 \\ 1.1456 \\ 4.85 \\ 4.85 \end{bmatrix}$$

In step 3, solving the regional upper-level optimization problem Eq.\ref{eq:Iter1_Step3} in the $CR_1$ produces the next iterate $x_2 = [-4.85, 1.2436]$.

\begin{equation}
\label{eq:Iter1_Step3}
\begin{aligned}
\min_{x_2} \quad & F(x_2) = \begin{bmatrix} 0 & 28.1 \end{bmatrix} x_2 + \begin{bmatrix} 0 & 26.2 \end{bmatrix} z \\
\text{s.t.} \quad & \begin{bmatrix} 1 & 0 \\ 0 & 1 \\ -1 & 0 \\ 0 & -1 \end{bmatrix} x_2 \le \begin{bmatrix} 4.85 \\ 4.85 \\ 4.85 \\ 4.85 \end{bmatrix} \\
&E_1x_2 \leq f_1\\
& z^1=0\\
&z^2=K_{C1}x_2+h_{C1}\\
&\begin{bmatrix} -2.0 & 5.9 \\ -3.0 & -3.0 \\ 7.1 & 4.7 \\ 0 & 2.0 \\ 0 & -10.0 \end{bmatrix} z \le \begin{bmatrix} 3 \\ 14.6 \\ 6.2 \\ 10.7 \\ 100 \end{bmatrix} + \begin{bmatrix} 0 & -7.8 \\ 0 & -9 \\ -4.7 & -2.4 \\ -7.8 & -5.9 \\ 0 & 0 \end{bmatrix} x_2 \\
\text{where} \quad & x_2 \in \mathbb{R}^2
\end{aligned}
\end{equation}

\subsubsection{Iteration 2}
In step 1, we fixed the $x_2$ and solved the lower-level problem, as shown in Eq.\ref{eq:Iter2_Step1}, which yields the optimal response $y_2 = [1, -0.7967]$ and a better upper objective value of $13.97$. Therefore, the $(x_2, y_2)$ is stored as the incumbent best solution $(x^*, y^*)$.

\begin{equation}
\label{eq:Iter2_Step1}
\begin{aligned}
\min_{y_2} \quad & f(y_2) = \begin{bmatrix} 1 & -3.6 \end{bmatrix} y_2\\
\text{s.t.} \quad & \begin{bmatrix} -2.0 & 5.9 \\ -3.0 & -3.0 \\ 7.1 & 4.7 \\ 0 & 2.0 \\ 0 & -10.0 \end{bmatrix} y_2 \le \begin{bmatrix} 3 \\ 14.6 \\ 6.2 \\ 10.7 \\ 100 \end{bmatrix} + \begin{bmatrix} 0 & -7.8 \\ 0 & -9 \\ -4.7 & -2.4 \\ -7.8 & -5.9 \\ 0 & 0 \end{bmatrix} x_2\\
\text{where} \quad & y^1_2 \in \{ 0, 1 \}, \quad y^2_2 \in \mathbb{R}
\end{aligned}
\end{equation}

In Step 2, sensitivity analysis identifies a new active set $\mathcal{A}_2 = \{0\}$, with the binary variable fixed at $y^1 = 1$. This defines the second critical region, $CR_2$, where the optimal solution is characterized by the linear mapping $y^2(x) = K_{C2}x + h_{C2}$. Where:
$$K_{C2} = \begin{bmatrix} 0 & -1.3220 \end{bmatrix}, \quad h_{C2} = \begin{bmatrix} 0.8475 \end{bmatrix}$$
The boundaries of this region are explicitly defined by the polyhedral constraints $E_2x \leq f_2$, where:

$$E_2 = \begin{bmatrix} 
0 & 1 \\ 
0.7765 & -0.6301 \\ 
-1 & 0 
\end{bmatrix}, \quad 
f_2 = \begin{bmatrix} 
1.5535 \\ 
-0.8068 \\ 
4.85
\end{bmatrix}
$$

In step 3, solving the regional upper-level optimization problem in the $CR_2$ as shown in Eq. \ref{eq:Iter2_Step3} produces the next iterate $x_3 = [-4.85, 1.5535]$.

\begin{equation}
\label{eq:Iter2_Step3}
\begin{aligned}
\min_{x_3} \quad & F(x_3) = \begin{bmatrix} 0 & 28.1 \end{bmatrix} x_3 + \begin{bmatrix} 0 & 26.2 \end{bmatrix} z \\
\text{s.t.} \quad & \begin{bmatrix} 1 & 0 \\ 0 & 1 \\ -1 & 0 \\ 0 & -1 \end{bmatrix} x_3 \le \begin{bmatrix} 4.85 \\ 4.85 \\ 4.85 \\ 4.85 \end{bmatrix} \\
&E_2x_3 \leq f_2\\
&z^1=1\\
&z^2=K_{C2}x_3+h_{C2}\\
&\begin{bmatrix} -2.0 & 5.9 \\ -3.0 & -3.0 \\ 7.1 & 4.7 \\ 0 & 2.0 \\ 0 & -10.0 \end{bmatrix} z \le \begin{bmatrix} 3 \\ 14.6 \\ 6.2 \\ 10.7 \\ 100 \end{bmatrix} + \begin{bmatrix} 0 & -7.8 \\ 0 & -9 \\ -4.7 & -2.4 \\ -7.8 & -5.9 \\ 0 & 0 \end{bmatrix} x_3 \\
\text{where} \quad & x_3 \in \mathbb{R}^2\\
\end{aligned}
\end{equation}

\subsubsection{Iteration 3}
In step 1, we fixed $x_3$ and solved the lower-level optimization problem as shown in Eq. \ref{eq:Iter3_Step1}, which yields the optimal response $y_3 = [1,-1.2063]$. Since the upper objective value ($12.0479$) is smaller than the incumbent best solution, the best solution is updated. 

\begin{equation}
\label{eq:Iter3_Step1}
\begin{aligned}
\min_{y_3} \quad & f(y_3) = \begin{bmatrix} 1 & -3.6 \end{bmatrix} y_3\\
\text{s.t.} \quad & \begin{bmatrix} -2.0 & 5.9 \\ -3.0 & -3.0 \\ 7.1 & 4.7 \\ 0 & 2.0 \\ 0 & -10.0 \end{bmatrix} y_3 \le \begin{bmatrix} 3 \\ 14.6 \\ 6.2 \\ 10.7 \\ 100 \end{bmatrix} + \begin{bmatrix} 0 & -7.8 \\ 0 & -9 \\ -4.7 & -2.4 \\ -7.8 & -5.9 \\ 0 & 0 \end{bmatrix} x_3\\
\text{where} \quad & y^1_3 \in \{ 0, 1 \}, \quad y^2_3 \in \mathbb{R}
\end{aligned}
\end{equation}

We check that $x_3,y_3^1$ is in the previously explored critical region $CR_2$. The algorithm converges, and report the best solution as $x_3 = [-4.85, 1.5535]$; $y_3 = [1,-1.2063]$; upper objective value $=12.0479$.

\section{Computational Implementation}
\label{Sec:Computational}

To evaluate the performance of the proposed Parametric Region Search (PRS) primal heuristic, a comparative analysis on time and solution quality was conducted against the Data-driven Optimization of bi-level Mixed-Integer NOnlinear (DOMINO) meta-heuristic solving framework \citep{beykal2020domino}. The DOMINO framework treats the lower-level problem as a black box and uses a metaheuristic approach to identify optimal upper-level variables that yield the best objective function value. Specifically, two derivative-free metaheuristic methods were utilized: the Constrained Optimization BY Linear Approximations (COBYLA) \cite{powell1994direct}, and Improved Stochastic Ranking Evolution Strategy (ISRES) \citep{runarsson2005search}. COBYLA functions as a local solver, constructed upon iterative linear approximations within a dynamically shrinking trust region. In contrast, ISRES is an evolution-based global optimization solver. The COBYLA and ISRES algorithms were implemented using NLopt 2.7.1 \citep{Johnson:NLopt}. The generation of critical regions required by the PRS heuristic is implemented utilizing PPOPT 1.6.0 \citep{kenefake2022ppopt}. All single-level MILP optimization subproblems were solved using Gurobi 13.0 \citep{gurobi2025}. 

Computational experiments were conducted on a Dell Pro Tower QCT1250 workstation equipped with an Intel Core Ultra 7 265 processor, featuring 20 cores and a base frequency of 2.40 GHz. The system is supported by 32.0 GB of installed physical memory (RAM) and operates on a 64-bit Windows 11 Education environment.

\subsection{Time comparison}
\label{sec:time}

\begin{figure}[htbp]
    \centering
    \includegraphics[width=0.75\linewidth]{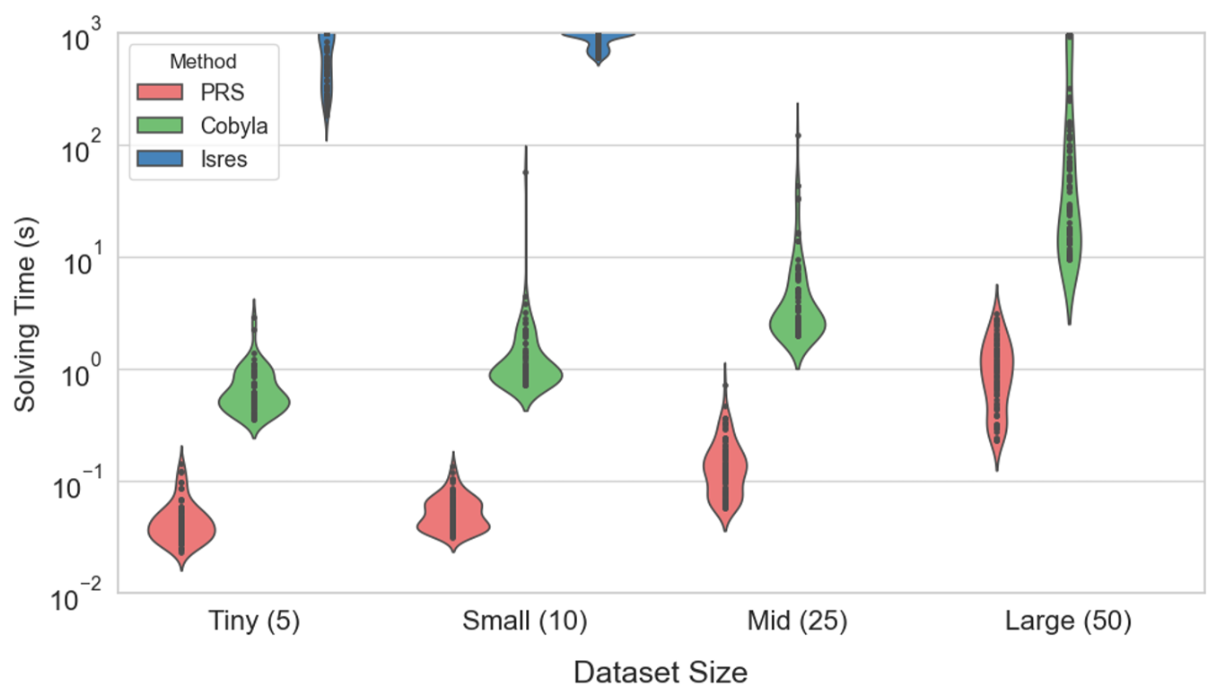}
    \caption{Algorithms time performance on comparison of different sizes: Tiny (5), Small (10), Mid (25), Large (50). Detailed size information is in Table \ref{tab:bilevel_sizes}.}
    \label{fig:Time}
\end{figure}

The solution time performance of the PRS, COBYLA and ISRES was assessed on $100$ randomly generated, non-trivial \citep{tsai2025assessing}, and feasible bilevel Mixed-Integer Linear Programming instances. The HPR response obtained by Algorithm \ref{alg:hprr} is given to all three algorithms as an initial point. To evaluate method scalability, four distinct problem structure sizes with a parameter density of 0.7 were utilized, as shown in Table \ref{tab:bilevel_sizes}. The problems are LP-MILP bilevel problems since the metaheuristic is tailored for continuous variables. The resulting solving times, capped at $1000$ seconds, are visualized using a violin plot (Figure~\ref{fig:Time}). The PRS was terminated according to the stopping criterion detailed in Section~\ref{sec:Parametric Region Search Heuristic}, whereas the DOMINO methods (COBYLA and ISRES) were executed until a relative error $10^{-4}$ was reached when solving lower-level optimization subproblem in Gurobi. Note that ISRES is not evaluated in mid and large cases, since it could not converge within the time budget.

\begin{table}[ht]
    \footnotesize
    \centering
    \caption{Problem Sizes for LP-MILP Bilevel Optimization Test Cases}
    \label{tab:bilevel_sizes}
    \begin{tabular}{|c|c|c|c|c|}
        \hline
        \textbf{Size (Total Number of Variables)} & \textbf{Tiny (5)} & \textbf{Small (10)} & \textbf{Mid (20)} & \textbf{Large (50)}\\
        \hline
        
        \multicolumn{5}{|l|}{\textbf{Upper-Level Problem (min $F(x,y)$)}} \\
        \hline
        Total Upper-Level Variables (Mixed) & 5 & 10 & 20 & 50 \\
        Continuous Upper-Level Variables ($|x_C|$) & 5 & 10 & 20 & 50 \\
        Binary Upper-Level Variables ($|x_I|$) & 0 & 0 & 0 & 0 \\
        Upper-Level Constraints (Bounds) & 10 & 20 & 40 & 100 \\
        \hline
        
        \multicolumn{5}{|l|}{\textbf{Lower-Level Problem ($y \in \arg\min f(x,z)$)}} \\
        \hline
        Total Lower-Level Variables (Mixed) & 5 & 10 & 20 & 50 \\
        Binary Lower-Level Variables ($|z_I|$) & 2 & 5 & 10 & 25 \\
        Continuous Lower-Level Variables ($|z_C|$) & 3 & 5 & 10 & 25 \\
        Lower-Level Constraints & 3 & 3 & 5 & 10 \\
        \hline
    \end{tabular}
    \vspace{0.5em}
\end{table}

Figure~\ref{fig:Time} illustrates the solving time distribution for the PRS method alongside the meta-heuristic optimizers COBYLA and ISRES, across the four problem scales. For all dataset sizes, the solution time distribution for PRS (red) is markedly lower than that for COBYLA (green) and ISRES (blue). The PRS distribution is approximately one order of magnitude lower than COBYLA and three orders of magnitude lower than ISRES, indicating rapid convergence within a small number of iterations. The ISRES method, a global solver, exhibits significantly higher solving times across the tested scales, with its distribution compressed near the maximum time cap of $10^3$ seconds. This comparative analysis demonstrates the {PRS} method's ability to generate a primal solution with superior computational efficiency across the tested scalability range. 

\begin{figure}[htbp]
    \centering
    \includegraphics[width=0.75\linewidth]{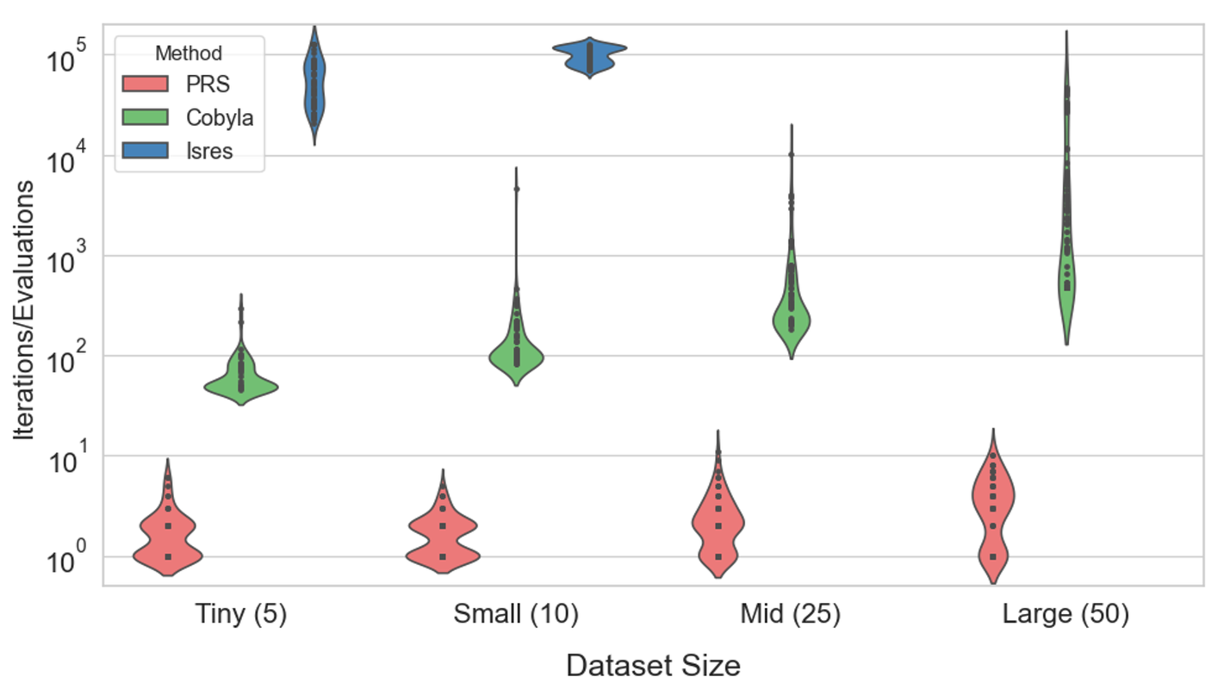}
    \caption{Algorithms comparison on number of iterations(PRS), and evaluations(COBYLA/ISRES) for different sizes: Tiny (5), Small (10), Mid (25), and Large (50). Detailed size information is in Table \ref{tab:bilevel_sizes}.}
    \label{fig:Evaluations}
\end{figure}

Figure~\ref{fig:Evaluations} illustrates the distribution of iteration counts for the PRS alongside the evaluation totals for COBYLA and ISRES. The PRS algorithm is achieving convergence within 10 iterations, whereas COBYLA and ISRES require hundreds to thousands of evaluations. Each PRS iteration requires solving three distinct optimization tasks: a lower-level Mixed-Integer Linear Program (MILP), a Linear Program (LP) for critical region generation, and a regional upper-level LP. While this involves only three subproblems, the total computational time per iteration typically exceeds that of solving three standard lower-level MILPs in COBYLA or ISRES. This performance gap occurs because PRS executes a sequence of heterogeneous problems. In contrast, COBYLA and ISRES solve a series of closely related MILP instances, which allows Gurobi to significantly reduce solve times through warm-starting.

\subsection{Quality comparison}
\label{Sec: Quality}
The nature of meta-heuristic solvers can generate better and better solutions before it converges by increasing solution time. However, it is not practical to let the heuristic run to converge for these solvers because the converging time of these methods is large. Therefore, people often implement a time cut to force terminate the meta-heuristic. In the following section we compared the solution quality of the PRS and the COBYLA/ISRES by changing different time cuts on small and large random bilevel problems in the Section \ref{sec:time}.

\begin{figure}[htbp]
    \centering
    \includegraphics[width=0.8\linewidth]{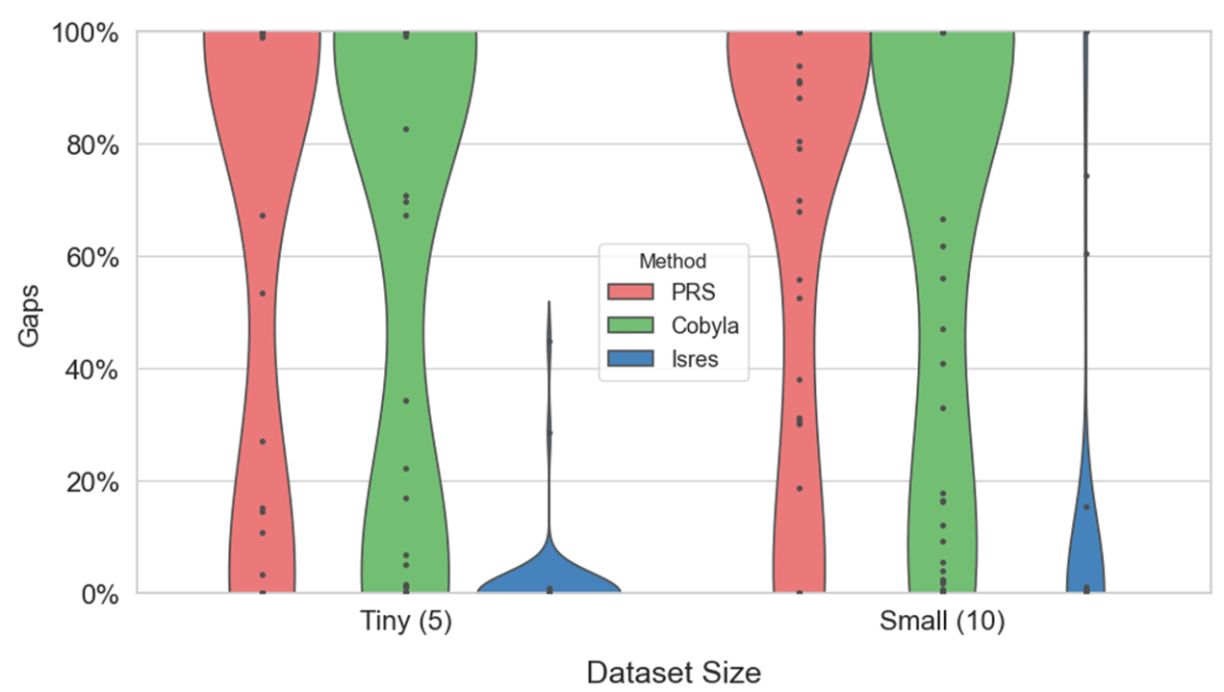}
    \caption{Algorithms gaps performance on comparison of different sizes: Tiny (5) and Small (10). Detailed size information is in Table \ref{tab:bilevel_sizes}.}
    \label{fig:Gaps}
\end{figure}

Figure \ref{fig:Gaps} is a violin plot that visualizes the distribution and density of optimality gaps across different dataset sizes. The ``width" of the violin represents the kernel density estimation, showing where the data points are most concentrated. The optimality gap is defined as:

\begin{equation}
Gap = 1 - \frac{F_{Alg}-F_{HPRR}}{F_{best}-F_{HPRR}}
\end{equation}
In this expression, $F_{Alg}$ denotes the objective value obtained by the specific algorithm being evaluated, while $F_{best}$ represents the best objective value found across all tested methods for a given instance. The definition of the gap assumes that at least one of the three evaluated methods attains a global optimum. This assumption is justified by the use of the global solver ISRES, as the best-obtained solutions are expected to be near-global optima for the Tiny and Small-size problem sets. Instances in which none of the three algorithms improved upon the initial High-Point Relaxation response (HPRR) are excluded from this analysis. ISRES achieves the lowest gaps, averaging $1\%$ and $9\%$ for the Tiny and Small-size sets, respectively. Notably, instances where ISRES exhibits larger gaps are primarily due to the gap definition; the absolute gap values remain small relative to the magnitude of the objective function. It is observed that in approximately half of the instances, both PRS and COBYLA yield an improvement over the baseline (HPR response). Therefore, PRS and COBYLA exhibit much higher gaps, averaging $64\%$ and $62\%$ for Tiny-size problems, and $70\%$ and $71\%$ for Small-size problems. The wide upper section of the violins for PRS and COBYLA indicates that in approximately half of the test instances, these methods failed to improve upon the baseline. The distribution of the results indicates that these bilevel problems remain challenging for the tested heuristics. The algorithms can sometimes stagnate at local optima, which prevents them from successfully converging toward lower objective values.

\begin{figure}[htbp]
    \centering
    \includegraphics[width=1\linewidth]{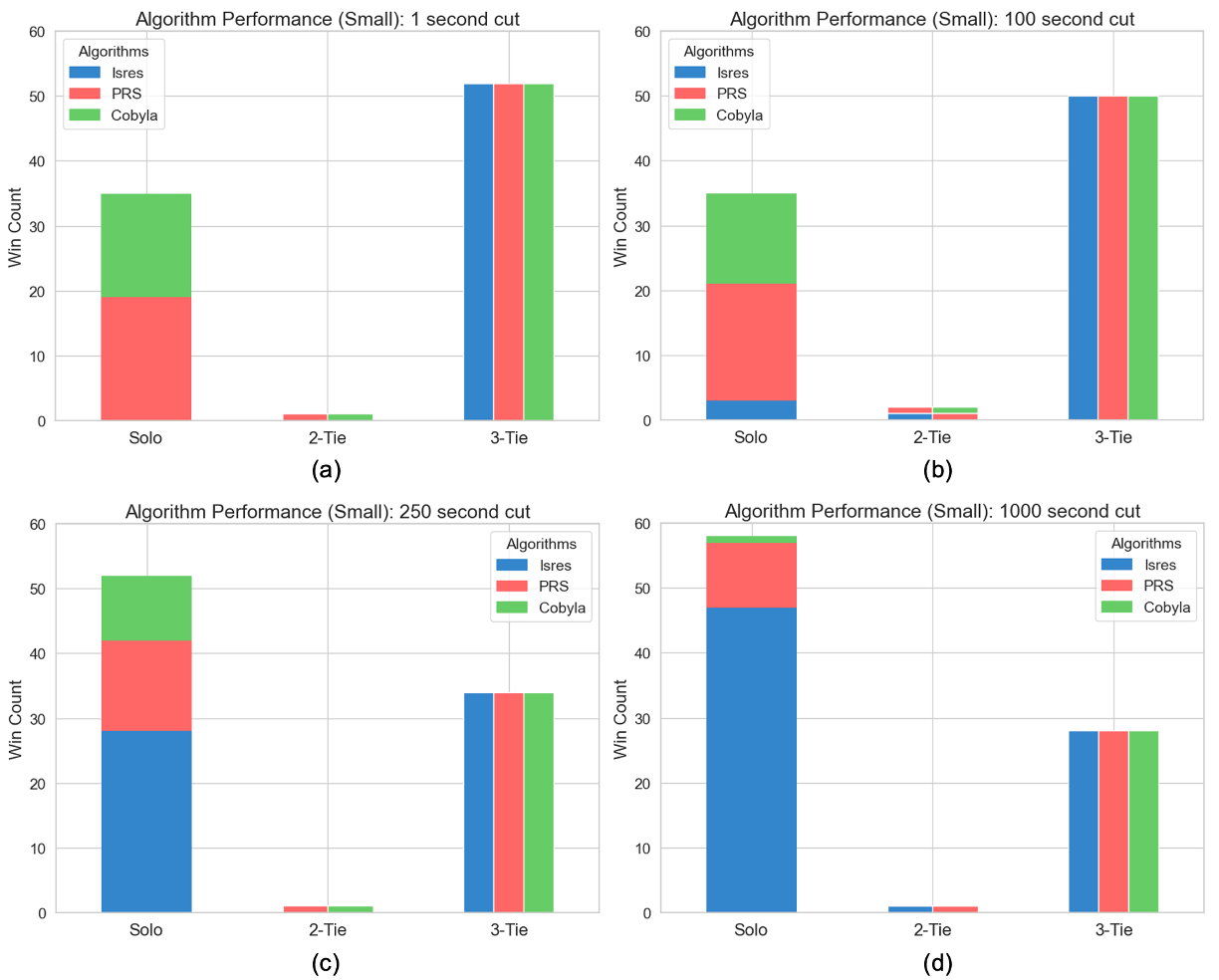}
    \caption{Algorithm performance comparison on small problems of different time cuts: (a) 1 second, (b) 100 second, (c) 250 second, (d) 1000 second.}
    \label{fig:Time cut}
\end{figure}

Figure \ref{fig:Time cut} illustrates the comparative solution quality of PRS, COBYLA, and ISRES across 100 randomly generated small-scale instances, where a win denotes the identification of a superior feasible objective within a fixed temporal budget. The results are presented by solo wins and multi-algorithm ties (more than one algorithm reaches the same optimal value) across four time horizons: 1s, 100s, 250s, and 1000s. Under tight constraints (1s), ISRES yields no solo wins, while PRS and COBYLA perform competitively; however, as the budget extends to 1000s, ISRES leverages its global search characteristics to secure nearly 50 solo wins, displacing the local solvers. Notably, the number of three-way ties diminishes from over 50 to under 30 as time increases, indicating that extended budgets allow for the ISRES to solve globally. The proposed PRS heuristic effectively utilizes parametric region structures to identify high-quality bilevel feasible solutions that evade COBYLA and, in specific instances, outperform ISRES even under a substantial computational budget.

\begin{figure}[htbp]
    \centering
    \includegraphics[width=1\linewidth]{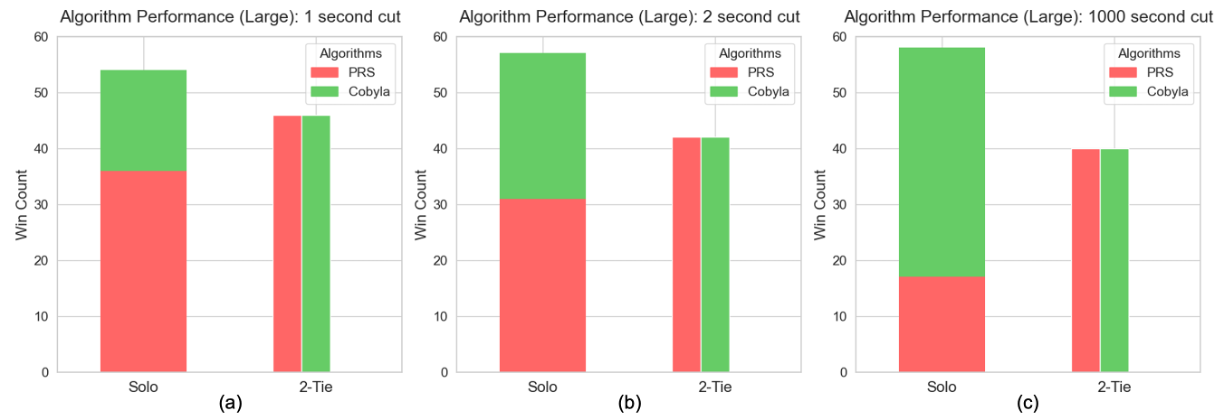}
    \caption{Algorithm performance comparison on large problems of different time cuts: (a) 1 second, (b) 2 second, (d) 1000 second.}
    \label{fig:Time cut2}
\end{figure}

Figure \ref{fig:Time cut2} illustrates the performance of PRS compared to COBYLA across 100 large-scale random instances, with COBYLA evaluated under computational budgets of 1, 2, and 1000 seconds. The ISRES algorithm is omitted from this comparison due to its high computational cost and failure to converge within the time constraints for high-dimensional problems. The results reveal a temporal performance shift: as the time budget increases, COBYLA's solo-win frequency rises to approximately 40 instances at the 1000-second mark, while the solo-win count for PRS diminishes to fewer than 20. This trend, coupled with the reduction in 2-Tie results, indicates that COBYLA identifies superior solutions when given larger time budget. Nevertheless, the findings demonstrate that PRS maintains a distinct competitive advantage by occasionally identifying feasible solutions that COBYLA fails to reach locally, even with extended execution times.

Although metaheuristics generally provide superior solutions given sufficient time, the PRS heuristic remains highly competitive by uncovering high-quality feasible points that standard local and global solvers overlook, especially within restricted temporal budgets.

\section{Enhanced bilevel optimization through PRS}

Having a good primal heuristic for bilevel optimization can be useful in finding near-optimal solution, when the actual optimal solution for bilevel optimization can be extremely difficult to find for large MILP case.

 \subsection{Combined PRS with COBYLA}

\begin{figure}[htbp]
    \centering
    \includegraphics[width=0.62\linewidth]{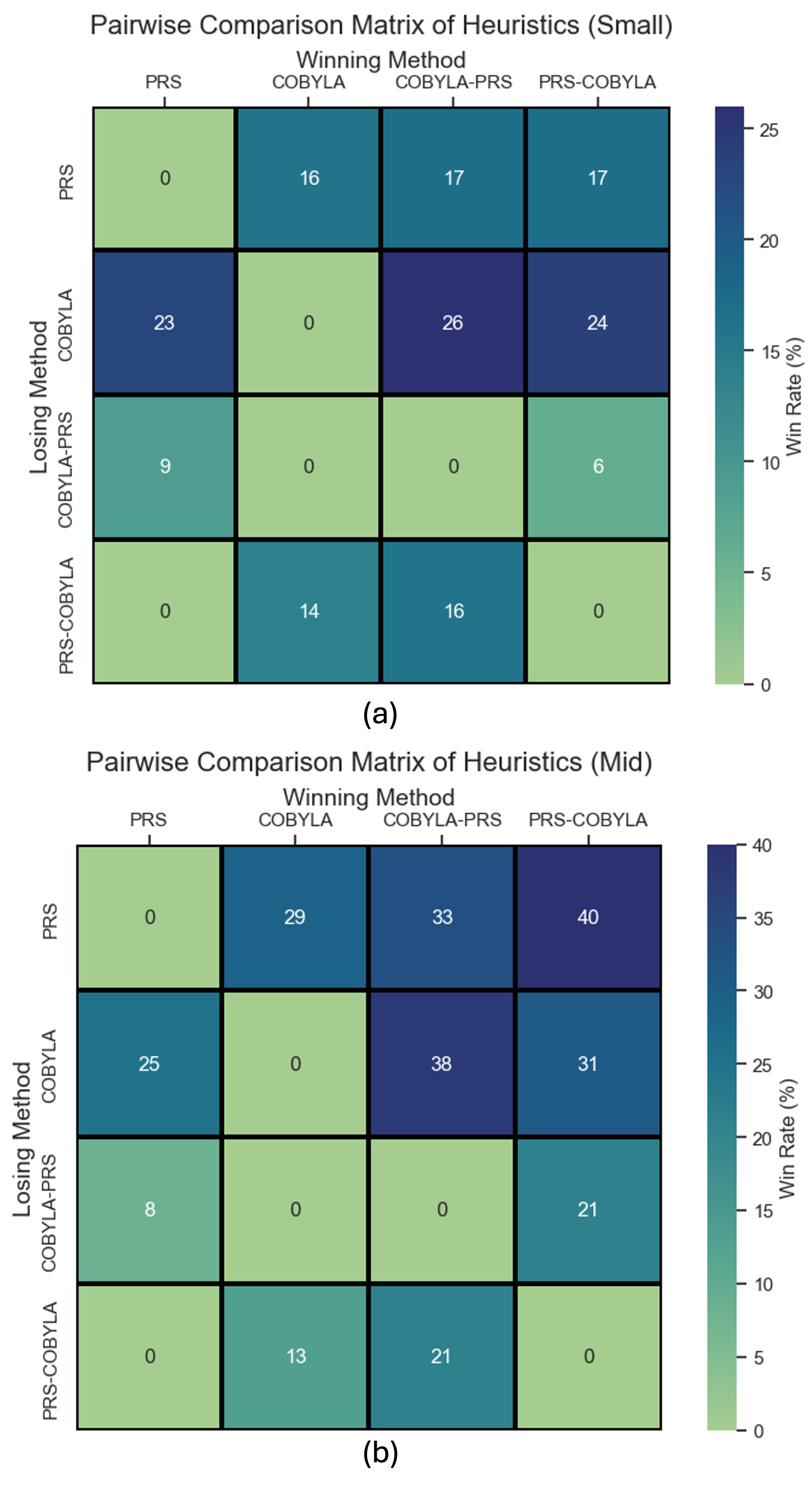}
    \caption{Pairwise win rate analysis of algorithm performance. The plots contrast the PRS, COBYLA, PRS-COBYLA, and COBYLA-PRS heuristics under two different sizes: (a) small-size and (b) mid-size problems.}
    \label{fig:pairwise-compare}
\end{figure}

This section shows the integration of the proposed PRS algorithms with the local derivative-free metaheuristic solver, COBYLA. By combining two heuristics that utilize distinct underlying logic to identify primal solutions, we aim to demonstrate a synergistic effect that enhances overall performance. Two hybrid frameworks were developed: PRS-COBYLA and COBYLA-PRS. In the PRS-COBYLA configuration, the process begins with an initial bilevel feasible solution (HPR response) that is first refined by the PRS algorithm; the resulting point then serves as the starting feasible point for further optimization via COBYLA. Conversely, the COBYLA-PRS method reverses this sequence, employing COBYLA to improve the initial feasible point (HPR response) before utilizing the PRS algorithm to conduct the final optimization phase. The results of the comparison of the two algorithms along with the original PRS and COBYLA are compared in Figure~\ref{fig:pairwise-compare}.

The pairwise performance comparison presented in Figure~\ref{fig:pairwise-compare} utilizes the small and mid-sized randomly generated problem sets detailed in Section~\ref{Sec:Computational}. A win is defined as achieving an objective value lower than the competing method by a margin exceeding $\epsilon = 10^{-5}$, and the win rate is defined as the number of wins divided by the total number of instances. Note that the sum of the win rates for any given pair does not necessarily equal $100\%$, as the remaining percentage corresponds to instances where both methods yield equivalent solutions (ties), which are not explicitly depicted in the figure. The computational time of the hybrid approach approximately equals the sum of its constituent algorithms. PRS and COBYLA demonstrate competitive efficacy; specifically, PRS achieves win rates of $23\%$ and $25\%$, while COBYLA has $16\%$ and $29\%$ in two different sized problem sets. These results aligned with the observations in Section~\ref{Sec:Computational}, showing that both algorithms effectively navigate the feasible region through different search mechanisms.

By construction, the performance of the COBYLA algorithm is bounded above by the COBYLA-PRS variant, and the PRS is bounded by the hybrid PRS-COBYLA. This relationship is a direct consequence of the sequential algorithmic architecture, wherein the secondary phase initiates from the best feasible solution identified during the primary phase, thereby guaranteeing that the hybrid result is at least as optimal as its precursor. The hybrid results shows this structural advantage: the COBYLA-PRS framework yields improvements over the standard COBYLA of $26\%$ and $38\%$ for small-size and medium-size instances, respectively. Similarly, the PRS-COBYLA configuration improves the solution quality of the PRS by $17\%$ for small-size problems and $40\%$ for medium-size problems.

While neither COBYLA-PRS nor PRS-COBYLA demonstrates strict dominance over their constituent non-hybrid counterparts, numerical experiments suggest that hybridization significantly increases the probability of attaining superior solutions. Specifically, COBYLA-PRS achieves win rates of $17\%$ and $33\%$ for small-size and medium-size instances, respectively, outperforming the $8\%$ and $9\%$ win rates observed for PRS. A similar trend is evident in PRS-COBYLA, which yields win rates of $24\%$ and $31\%$, whereas the standalone COBYLA algorithm achieves only $14\%$ and $13\%$ across the same categories. When comparing the two hybrid configurations, their performances are relatively competitive; PRS-COBYLA maintains win rates of $6\%$ and $21\%$, while COBYLA-PRS records $16\%$ and $21\%$ for small and medium instances. To sum up, the two hybrid strategies are comparable and offer advantages over vanilla PRS and COBYLA methods.

\section{Conclusion and Discussion}
This work proposes the Parametric Region Search (PRS), a novel primal heuristic designed to generate high-quality solutions for Mixed-Integer Bilevel Optimization (MIBO) at a low computational cost. Inspired by multiparametric programming, the algorithm leverages the multiparametric reformulation structure of MIBO problems to accelerate the search process.

Our study evaluates the solution quality and runtime of the PRS against two metaheuristic-based methods: DOMINO-COBYLA (local) and DOMINO-ISRES (global). According to the results, PRS can sometimes find better feasible solutions in a shorter amount of time and with fewer subproblems than both DOMINO variants. Furthermore, we developed two hybrid approaches by integrating PRS with DOMINO-COBYLA. Both hybrids (COBYLA-PRS and PRS-COBYLA) outperform their standalone counterparts in terms of solution quality, with a small add-on time.

The PRS algorithm is a versatile framework with significant potential for bilevel optimization research. Future work will focus on adapting the method to generate high-quality upper bounds for other solvers and scaling the approach to address larger and more complex problems.

\section*{Declaration of interest}
None
\section*{Declaration of generative AI and AI-assisted technologies in the writing process.}
During the preparation of this work, the authors used Gemini in order to improve the readability and language of the manuscript. After using this tool/service, the authors reviewed and edited the content as needed and took full responsibility for the content of the published article.

\section*{Acknowledgments}
This material is based upon work supported by the National Science Foundation under NSF Award \#2328160 and the Department of Chemical and Biological Engineering at the University of Wisconsin-Madison. 

{\small\bibsep=0pt
\bibliography{example}
}

\end{document}